\pdfoutput=1
\RequirePackage{ifpdf}
\ifpdf 
\documentclass[pdftex]{sigma}
\else
\documentclass{sigma}
\fi

\numberwithin{equation}{section}

\usepackage{tikz-cd}
\usepackage{mathtools}

\newtheorem{tw}{Theorem}[section]
\newtheorem*{Theorem*}{Theorem}
\newtheorem{cor}[tw]{Corollary}

\newtheorem{prop}[tw]{Proposition}

\theoremstyle{definition}
\newtheorem{defi}[tw]{Definition}

\newtheorem{egg}[tw]{Example}
\newtheorem{rem}[tw]{Remark}
\newtheorem{Not}[tw]{Notation}
\newtheorem{ques}[tw]{Question}

\begin{document}

\allowdisplaybreaks

\newcommand{\arXivNumber}{2506.12422}

\renewcommand{\PaperNumber}{072}

\FirstPageHeading

\ShortArticleName{The Spectral Sequence for Actions of Abelian Lie Groups}

\ArticleName{A Study of the Spectral Sequence for Locally Free \\ Isometric Actions of Abelian Lie Groups}

\Author{Pawe{\l} RA\'ZNY}

\AuthorNameForHeading{P.~Ra\'zny}

\Address{Institute of Mathematics, Jagiellonian University in Cracow,\\ ul. prof. Stanislawa {\L}ojasiewicza 6, 30-348 Krak\'ow, Poland}
\Email{\mail{pawel.razny@uj.edu.pl}}

\ArticleDates{Received June 17, 2025, in final form July 27, 2026; Published online August 08, 2026}

\Abstract{We give an upper bound on the number of the page on which the spectral sequence corresponding to a locally free isometric action of an abelian Lie group degenerates. We give examples showing that these bounds are indeed sharp. Finally, we further justify the study of this sequence by exhibiting a potential application to the study of harmonic forms.}

\Keywords{Lie group actions; Lie algebra actions; foliations; basic cohomology; isometries}

\Classification{53C12}

\section{Introduction}
In~\cite{My4}, we have studied a Serre like spectral sequence for locally free isometric Lie group actions on compact manifolds by showing the following result.

\begin{tw}\label{RecG} Let $(M^{n+s},g)$ be a compact manifold with an isometric locally free action of an~$s$-dimensional connected Lie group $G$. Then, there is a spectral sequence $E^{p,q}_r$ with
\begin{itemize}\itemsep=0pt
\item $E^{p,q}_2=H^p(M\slash\mathcal{F}, H^{q}(\mathfrak{g}))$, where $\mathfrak{g}$ is the Lie algebra of $G$ and $\mathcal{F}$ is the foliation generated by the fundamental vector fields of the action $\xi_1,\dots,\xi_s$.
\item $E^{p,q}_r$ converges to $H^{\bullet}_{\rm dR}(M)$.
\end{itemize}
\end{tw}
The purpose of this article is to study this spectral sequence in the special case when $G$ is abelian. In this case the sequence simplifies greatly (similarly as in the more restrictive case considered in~\cite{My3}) and this allows us to extract some additional information from it. Notably this simple case already has some interesting applications in differential geometry such as to the case of $\mathcal{K}$-structures (cf.\ \cite{B}) which we already explored in~\cite{My3} or the case of $q$-contact manifolds of~\cite{F} which provided a major motivation for this article. Our main result is an upper bound on the number of page at which the spectral sequence degenerates based on the basic cohomology classes of ${\rm d}\eta_i$ (where $\eta_i$ denote $1$-forms corresponding to the generators of $\mathfrak{g}^*$ as explained in the subsequent section).

\begin{tw}\label{Main} Let $(M^{n+s},g)$ be a compact Riemannian manifold with a locally free isometric action of a connected abelian Lie group $G$ $($of dimension $s)$ such that $\{{\rm d}\eta_1,\dots,{\rm d}\eta_s\}$ has basic cohomological rank $k$ $\bigl($i.e., $\dim \langle [{\rm d}\eta_1], \dots, [{\rm d}\eta_s] \rangle = k$ in $H^2(M\slash\mathcal{F})\bigr)$. Then, the spectral sequence~$E^{p,q}_r$ degenerates at the latest at the $(k+2)$-th page.
\end{tw}
 This result is designed to shorten computations using this sequence by reducing the number of pages to consider. We follow this result in the subsequent section by a number of examples showing that the estimate provided is indeed sharp.

 The penultimate section is designed to motivate the above computational results by providing some justification for the value of this sequence. In it we explore a possible use of this sequence in studying harmonic forms on manifolds with isometric abelian Lie group actions. For this purpose, we provide an explicit example of computation of harmonic forms using our sequence. We believe that the correspondence between harmonic forms and basic harmonic forms can be still strengthened for group actions with more rigid geometric structures (which is justified by the full description given in the case of $\mathcal{S}$-structures in~\cite{My3}). Further study of this theory in these cases (e.g., $\mathcal{K}$-structures, $q$-contact structures) should produce interesting results.

In the final section, we show that the exact same bound applies to the spectral sequence of the Riemannian foliation given by the orbits of the action. We are grateful to the anonymous referees for pointing out this significant consequence of our work.
\section{Preliminaries}
\subsection{Foliations}
We provide a quick review of transverse structures on foliations. Most of the content of this section can be found, e.g., in~\cite{E1} (see also~\cite{M} for more information).
\begin{defi} A codimension $q$ foliation $\mathcal{F}$ on a smooth $n$-manifold $M$ is given by the following~data:
\begin{itemize}\itemsep=0pt
\item an open cover $\mathcal{U}:=\{U_i\}_{i\in I}$ of M,
\item a $q$-dimensional smooth manifold $T_0$,
\item for each $U_i\in\mathcal{U}$, a submersion $f_i\colon U_i\rightarrow T_0$ with connected fibers (these fibers are called plaques),
\item for all intersections $U_i\cap U_j\neq\varnothing$, a local diffeomorphism $\gamma_{ij}$ of $T_0$ such that $f_j=\gamma_{ij}\circ f_i$.
\end{itemize}
The last condition ensures that plaques glue nicely to form a partition of $M$ consisting of submanifolds of $M$ of codimension $q$. This partition is called a foliation $\mathcal{F}$ of $M$ and the elements of this partition are called leaves of $\mathcal{F}$.
\end{defi}
We call \smash{$T=\coprod_{U_i\in\mathcal{U}}f_i(U_i)$} the transverse manifold of $\mathcal{F}$. The local diffeomorphisms $\gamma_{ij}$ generate a pseudogroup $\Gamma$ of transformations on $T$ (called the holonomy pseudogroup). The space of leaves $M\slash\mathcal{F}$ of the foliation $\mathcal{F}$ can be identified with $T\slash\Gamma$.
\begin{defi}
 A smooth form $\omega$ on $M$ is called transverse if for any vector field $X\in\Gamma (T\mathcal{F})$ (where $T\mathcal{F}$ denotes the bundle tangent to the leaves of $\mathcal{F}$) it satisfies $\iota_X\omega=0$. Moreover, if additionally ${\rm d}\omega=0$ is also transverse, then $\omega$ is said to be basic. Basic $0$-forms will be called basic functions henceforth.
\end{defi}
Basic forms are in one to one correspondence with $\Gamma$-invariant smooth forms on $T$. It is clear that ${\rm d}\omega$ is basic for any basic form $\omega$. Hence, the set of basic forms of $\mathcal{F}$ (denoted $\Omega^{\bullet}(M\slash\mathcal{F})$) is a subcomplex of the de Rham complex of $M$. We define the basic cohomology of $\mathcal{F}$ to be the cohomology of this subcomplex and denote it by $H^{\bullet}(M\slash\mathcal{F})$. A transverse structure to $\mathcal{F}$ is a~$\Gamma$-invariant structure on $T$. Among such structures the following is most relevant to our work.

\begin{defi}
A~foliation $\mathcal{F}$ is said to be Riemannian if~$T$ has a $\Gamma$-invariant Riemannian metric. This is equivalent to the existence of a Riemannian metric~$g$ (called the transverse Riemannian metric) on $N\mathcal{F}:=TM\slash T\mathcal{F}$ with $\mathcal{L}_Xg=0$ for all vector fields~$X$ tangent to the leaves.\looseness=-1
\end{defi}
\begin{defi}
A~foliation $\mathcal{F}$ is said to be transversely orientable if $T$ is orientable and all the~$\gamma_{ij}$ are orientation preserving. This is equivalent to the orientability of~$N\mathcal{F}$.
\end{defi}
These structures enable the construction of a transverse version of Hodge theory (see~\cite{E1}). Firstly, we recall a special class of Riemannian foliations on which the aforementioned theory is greatly simplified.

\begin{defi} A codimension $q$ foliation $\mathcal{F}$ on a compact connected manifold $M$ is called homologically orientable if $H^{q}(M\slash\mathcal{F})=\mathbb{R}$. A foliation $\mathcal{F}$ on a compact manifold $M$ is called homologically orientable if its restriction to each connected component of $M$ is homologically orientable.
\end{defi}
We will see at the end of this section that foliations considered in this paper are homologically orientable (under the assumption that $M$ is orientable) and hence we shall restrict our exposition to this case.

Let us firstly note that for Riemannian foliations (on compact manifolds) homological orientability implies transverse orientability (the converse is not true see, e.g.,~\cite{E1}). We believe this is known and it is implied in~\cite{E1} but since we could not find a proof we elect to include one for completness.
\begin{prop} Let $\mathcal{F}$ be a homologically orientable Riemannian foliation on a compact manifold $M$. Then $\mathcal{F}$ is transversely orientable.
\end{prop}

\begin{proof} Without loss of generality, we can assume that $M$ is connected and lets assume to the contrary that $\mathcal{F}$ is not transversely orientable. This is equivalent to $V:=\bigwedge^q (N\mathcal{F})^*$ being non-trivial. Consider the double cover $p\colon \tilde{M}\to M$ given by the total space of the sphere bundle of~$V$. The foliation $\mathcal{F}$ lifts to a foliation $\tilde{\mathcal{F}}$ on $\tilde{M}$ along with the transverse Riemannian metric. Moreover, $\tilde{\mathcal{F}}$ is transversely orientable (since \smash{$\bigwedge^q \bigl(N\tilde{\mathcal{F}}\bigr)^*=p^*V$} is trivial) and homologically orientable. To see the latter, note that since \smash{$\Omega^k(M\slash\mathcal{F})\cong\Omega^k\bigl(\tilde{M}\slash\tilde{\mathcal{F}}\bigr)^{\mathbb{Z}_2}$} and taking the invariant part with respect to a~finite group action commutes with taking cohomology we can see that \smash{$\mathbb{R}\cong H^q(M\slash\mathcal{F})\cong H^q\bigl(\tilde{M}\slash\tilde{\mathcal{F}}\bigr)^{\mathbb{Z}_2}\subset H^q\bigl(\tilde{M}\slash\tilde{\mathcal{F}}\bigr)$}. By the well-known fact (see~\cite{E1}) that the top basic cohomology of a Riemannian foliation on a compact connected manifold is either~$0$ or~$\mathbb{R}$, we conclude that \smash{$H^q\bigl(\tilde{M}\slash\tilde{\mathcal{F}}\bigr)^{\mathbb{Z}_2}= H^q\bigl(\tilde{M}\slash\tilde{\mathcal{F}}\bigr)\cong\mathbb{R}$} as declared earlier. Using the basic version of Poincar\'e duality (see~\cite{E1}) which we later recall as part of Theorem~\ref{HD}, we get that \smash{$H^q\bigl(\tilde{M}\slash\tilde{\mathcal{F}}\bigr)$} is generated by a basic nowhere vanishing form (which can be taken as the basic Hodge star of $1$). Averaging this form over the $\mathbb{Z}_2$ action and using \smash{$H^q\bigl(\tilde{M}\slash\tilde{\mathcal{F}}\bigr)^{\mathbb{Z}_2}= H^q\bigl(\tilde{M}\slash\tilde{\mathcal{F}}\bigr)$} (which guarantees that at no point the $\mathbb{Z}_2$ action takes this form to its opposite), we get the desired nowhere vanishing section of~$V$.
\end{proof}

Let $\mathcal{F}$ be a homologically orientable Riemannian foliation on a compact manifold $M$. One can use the transverse Riemannian metric along with a choice of orientation of $N\mathcal{F}$ to define the basic Hodge star operator \smash{$*_b\colon \Omega^k(M\slash\mathcal{F})\to\Omega^{q-k}(M\slash\mathcal{F})$} pointwise. This in turn allows us to define the basic adjoint operator
\[\delta_b=(-1)^{q(k+1)+1}*_b{\rm d}*_b.\]
\begin{rem} While we choose this to be the definition of $\delta_b$, it is in fact an adjoint of ${\rm d}$ with respect to an appropriate inner product on forms induced by the transverse metric $g$. However, the definition of this inner product is quite involved and not necessary for our purpose. Although, we shall state some of the classical results of basic Hodge theory which use this inner product. For details see~\cite{E1}.
\end{rem}
Using $\delta_b$, we can define the basic Laplace operator via
\[\Delta_b={\rm d}\delta_b+\delta_b{\rm d}.\]
As it turns out this operator has some nice properties similar to that of the classical Laplace operator. In particular, it is transversely elliptic in the following sense.

\begin{defi}
A basic differential operator of order $m$ is a linear map ${D\colon\! \Omega^{\bullet}(M\slash\!\mathcal{F})\!\rightarrow\!\Omega^{\bullet}(M\slash\!\mathcal{F})}$ such that in local coordinates $(x_1,\dots,x_p,y_1,\dots,y_q)$ (where $x_i$ are leaf-wise coordinates and~$y_j$ are transverse ones) it has the form
\begin{equation*}
D=\sum\limits_{|s|\leq m}a_s(y)\frac{\partial^{|s|}}{\partial^{s_1}y_1\cdots \partial^{s_q}y_q},
\end{equation*}
where $a_s$ are matrices of appropriate size with basic functions as coefficients. A basic differential operator is called transversely elliptic if its principal symbol is an isomorphism at all points of~${x\in M}$ and all non-zero, transverse, cotangent vectors at~$x$.
\end{defi}

In particular, this implies the following important result from~\cite{E1}.

\begin{tw}\label{HD} Let $\mathcal{F}$ be a Riemannian homologically orientable foliation on a compact manifold~$M$.~Then
\begin{enumerate}\itemsep=0pt
\item[$1.$] $H^{\bullet}(M\slash\mathcal{F})$ is isomorphic to the space of basic harmonic forms $\operatorname{Ker}(\Delta_b)$. In particular, it is finite-dimensional.
\item[$2.$] The basic Hodge star induces an isomorphism between $H^{k}(M\slash\mathcal{F})$ and $H^{q-k}(M\slash\mathcal{F})$ given by taking the class of the image through $*_b$ of a harmonic representative.
\item[$3.$] $($Hodge decomposition$)$ The space of basic forms splits orthogonally into
\[\Omega^{\bullet}(M\slash\mathcal{F})=\operatorname{Ker}(\Delta_b)\oplus \operatorname{Im}(\delta_b)\oplus \operatorname{Im}\bigl({\rm d}|_{\Omega^{\bullet}(M\slash\mathcal{F})}\bigr).\]
\end{enumerate}
\end{tw}
We finish this section by recalling the spectral sequence of a Riemannian foliation (cf.~\cite{ss1,ss2,ss3}).
\begin{defi}
We put
\[F^k_{\mathcal{F}}\Omega^r(M):=\{\alpha\in\Omega^r(M)\mid \iota_{X_{r-k+1}}\cdots \iota_{X_1}\alpha=0\text{ for } X_1,\dots,X_{r-k+1}\in\Gamma(T\mathcal{F})\}.\]
An element of \smash{$F^k_{\mathcal{F}}\Omega^r(M)$} is called an $r$-differential form of filtration $k$.
\end{defi}
The definition above in fact gives a filtration of the de Rham complex by differential ideals. Hence, via known theory from homological algebra we can construct a spectral sequence as~follows:
\begin{enumerate}\itemsep=0pt
\item The $0$-th page is given by \smash{$E_0^{p,q}=F_{\mathcal{F}}^{p}\Omega^{p+q}(M)\slash F_{\mathcal{F}}^{p+1}\Omega^{p+q}(M)$} and \smash{${\rm d}^{p,q}_0\colon E_0^{p,q}\to E^{p,q+1}_0$} is simply the morphism induced by ${\rm d}$.
\item The $r$-th page is given inductively by
\begin{align*}
E_{r}^{p,q}:={}&\operatorname{Ker}\bigl({\rm d}^{p,q}_{r-1}\bigr)\slash \operatorname{Im}\bigl({\rm d}^{p-r+1,q+r-2}_{r-1}\bigr)\\
={}&\frac{\bigl\{\alpha\in F^{p}_{\mathcal{F}}\Omega^{p+q}(M)\mid {\rm d}\alpha\in F_{\mathcal{F}}^{p+r}\Omega^{p+q+1}(M)\bigr\}}{F^{p+1}_{\mathcal{F}} \Omega^{p+q}(M)+{\rm d}\bigl(F_{\mathcal{F}}^{p-r+1}\Omega^{p+q-1}(M)\bigr)}.
\end{align*}
\item The $r$-th coboundary operator ${\rm d}_r\colon E_r^{p,q}\to E^{p+r,q-r+1}_r$ is again just the map induced by~${\rm d}$ (due to the description of the $r$-th page this has the target specified above and is well defined).
\end{enumerate}
Furthermore, since the filtration is bounded this spectral sequence converges and its final page is isomorphic to the cohomology of the cochain complex (in this case, the de Rham cohomology of $M$).
\begin{rem} The above spectral sequence can be thought of as a generalization of the Leray--Serre spectral sequence in de Rham cohomology to arbitrary Riemannian foliations (as opposed to fiber bundles).
\end{rem}

\subsection{The spectral sequence of a locally free action}
In this section we will recall the construction of the spectral sequence which is our main object of study. More precisely, we will recall the construction of the spectral sequence from Theorem~\ref{RecG} (from~\cite{My4}) in the special case when $G$ is abelian which is somewhat simpler (and follows closely the exposition in~\cite{My3}) and sufficient for the study conducted in this article. Let us start by recalling the following definition.

\begin{defi} We say that a Lie group action of $G$ on $M$ is locally free if all its isotropy groups are discrete.
\end{defi}

Let $M$ be a compact Riemannian manifold with a Riemannian metric $g$ and an isometric locally free action of an $s$-dimensional abelian Lie group G (without loss of generality, by taking the quotient, we can treat this group as a subgroup of $\operatorname{Diff}(M)$). Then the fundamental vector fields $\xi_1,\dots,\xi_s$ of this action span the involutive subbundle $T\mathcal{F}\subset TM$ (here by~$\mathcal{F}$ we denote the foliation on $M$ given by the orbits of the action) and satisfy $[\xi_i,\xi_j]=0$. This implies the following.

\begin{prop} Let $(M^{n+s},g)$ be a compact Riemannian manifold with a locally free action of a connected abelian Lie group $G\subset \operatorname{Diff}(M)$ by isometries. Then the closure $\overline{G}$ of $G$ in $\operatorname{Diff}(M)$ is a torus contained in the group $\operatorname{Isom}(M)$ of isometries on $M$.
\end{prop}

\begin{proof} Since the action of $G$ is isometric, we have $G\subset \operatorname{Isom}(M)$ which is known to be a finite-dimensional compact Lie group. Moreover, since $G$ is abelian, its closure is a compact abelian group and hence a torus.
\end{proof}

The next step is to classify forms on $M$ which are invariant under the action of $\overline{G}$. For this let us consider the $1$-forms $\eta_1,\dots,\eta_s$ defined by $\eta_i(\xi_j)=\delta_{ij}$ and $\eta_i\bigl(T\mathcal{F}^{\perp}\bigr)=0$. Despite the following result being proved almost verbatim as an analogous result in~\cite{My3}, we elect to recall it in full for the reader's convenience.

\begin{prop}\label{InvFor} Let $(M^{n+s},g)$ be a compact Riemannian manifold with a locally free action of a connected abelian Lie group $G\subset \operatorname{Diff}(M)$ by isometries. Then the forms ${\rm d}\eta_i$ are basic. Moreover, the following conditions are equivalent:
\begin{enumerate}\itemsep=0pt
\item[$(1)$] $\alpha$ is a $\overline{G}$-invariant form on $M$.
\item[$(2)$] \smash{$\alpha=\alpha_0+\sum_{k=1}^p\sum_{1\leq i_1<\cdots <i_k\leq s}\eta_{i_1}\cdots \eta_{i_k}\alpha_{i_1,\dots,i_k}$}, where $\alpha_0$ and $\alpha_{i_1,\dots,i_k}$ are basic for all indices $1\leq i_1<\cdots <i_k\leq s$.
\end{enumerate}
\end{prop}
\begin{proof}
Firstly, let us show that the forms ${\rm d}\eta_i$ are indeed basic. To see this, let us note that since $\overline{G}$ acts by isometries and the flows of $\xi_i$ preserve $T\mathcal{F}$, they have to preserve $T\mathcal{F}^{\perp}$ as well. Consider the formula
\[{\rm d}\eta_i(X_0,X_1)= X_0.(\eta_i(X_1))-X_1.(\eta_i(X_0)) -\eta_i([X_0,X_1]).\]
To prove that this form indeed vanishes on $T\mathcal{F}$, it suffices to consider $X_0$ to be an $\mathbb{R}$-linear combination of $\xi_1,\dots,\xi_p$ (by tensoriality). In this case, $\eta_i(X_0)$ is constant and hence the second term vanishes. On the other hand, we can without loss of generality (by tensoriality) take $X_1$ to be a sum of an $\mathbb{R}$-linear combination of $\xi_1,\dots,\xi_p$ and a vector field in $T\mathcal{F}^{\perp}$. Hence, similarly as before, the first term in the sum vanishes. Moreover, the final term also vanishes since $[X_0,X_1]$ has to be a section of $T\mathcal{F}^{\perp}$ (since the brackets of the form $[\xi_i,\xi_j]$ vanish and $T\mathcal{F}^{\perp}$ is preserved by the flows of $\xi_i$ as already remarked). Finally, for any $X\in\Gamma(T\mathcal{F})$, we have{\samepage
\[\mathcal{L}_X{\rm d}\eta_j={\rm d}\iota_X{\rm d}\eta_j+\iota_X{\rm d}{\rm d}\eta_j=0,\]
which implies that ${\rm d}\eta_j$ is indeed basic.}

Now assume that the second condition is true. Then it can be easily computed that for any~$\xi_j$ the equality $\mathcal{L}_{\xi_j}\alpha=0$ holds. Which in turn implies that $\alpha$ is $G$-invariant and consequently $\overline{G}$-invariant.

Now let us write the invariant form $\alpha$ as
\[\alpha=\alpha_0+\sum\limits_{k=1}^s\sum\limits_{1\leq i_1<\cdots <i_k\leq s}\eta_{i_1}\cdots \eta_{i_k}\alpha_{i_1,\dots,i_k},\]
where $\alpha_{i_1,\dots,i_k}$ are transverse for all indices $1\leq i_1<\cdots <i_k\leq s$. Due to the well-known formula
\[\mathcal{L}_{X}i_{Y}-i_{Y}\mathcal{L}_{X}=i_{[X,Y]},\]
we get that $i_{\xi_i}$ and $\mathcal{L}_{\xi_j}$ commute for $i,j\in\{1,\dots,s\}$. We shall now prove that the forms~$\alpha_0$ and~${\alpha_{i_1,\dots,i_k}}$ are basic by reverse induction on the number of indices. Hence, we start by proving that $\alpha_{1,\dots,s}$ is basic. Since $\alpha$ is basic and the vector fields $\xi_i$ are Killing, we have for any ${i\in\{1,\dots,s\}}$ the following equalities:
\[0=\mathcal{L}_{\xi_i}\alpha=i_{\xi_s}i_{\xi_{s-1}}\cdots i_{\xi_1}\mathcal{L}_{\xi_i}\alpha=\mathcal{L}_{\xi_i}i_{\xi_{s}}i_{\xi_{s-1}}\cdots i_{\xi_1}\alpha=\mathcal{L}_{\xi_i}\alpha_{1,\dots,s},\]
which proves that $\alpha_{1,\dots,s}$ is basic.

For the induction step, let us assume that all the $\alpha_{i_1,\dots,i_{k}}$ for $p\geq k>K$ are basic. We shall show that all $\alpha_{i_1,\dots,i_{K}}$ are basic as well. Using the assumption, we get for any $i\in\{1,\dots,s\}$ the following equalities:
\[0=\mathcal{L}_{\xi_i}\alpha=i_{\xi_{i_K}}i_{\xi_{i_{K-1}}}\cdots i_{\xi_{i_1}}\mathcal{L}_{\xi_i}\alpha=\mathcal{L}_{\xi_i}i_{\xi_{i_K}}i_{\xi_{i_{K-1}}}\cdots i_{\xi_{i_1}}\alpha=\mathcal{L}_{\xi_i}\alpha_{i_1,\dots,i_K},\]
which proves that $\alpha_{i_1,\dots,i_K}$ are basic for any set of indices $1\leq i_1<\cdots <i_k\leq s$.
\end{proof}

\begin{rem} Note that the induction assumption is used to pass to the final equality as it implies that all the terms with a greater number of indices than $K$ vanish under $\mathcal{L}_{\xi_i}$ as
\[\mathcal{L}_{\xi_i}\eta_{j_1}\cdots \eta_{j_k}\alpha_{j_1,\dots,j_k}=\eta_{j_1}\cdots \eta_{j_k}\mathcal{L}_{\xi_i}\alpha_{j_1,\dots,j_k}=0.\]
The first equality is due to the fact that $\mathcal{L}_{\xi_i}\eta_j=i_{\xi_i}{\rm d}\eta_j+{\rm d}(i_{\xi_i}\eta_j)=0$.
\end{rem}
Finally, we note that similarly as for the spectral sequence of a Riemannian foliation we have a filtration of the cochain complex of invariant forms $\Omega^r_{\overline{G}}(M)$ given by
\[F^k_{\mathcal{F}}\Omega^r_{\overline{G}}(M):=F^k_{\mathcal{F}}\Omega^r(M)\cap\Omega^r_{\overline{G}}(M).\]
Hence, via known theory from homological algebra, we can construct a spectral sequence as follows:
\begin{enumerate}\itemsep=0pt
\item The $0$-th page is given by \smash{$E_0^{p,q}=F_{\mathcal{F}}^{p}\Omega^{p+q}_{\overline{G}}(M)\slash F_{\mathcal{F}}^{p+1}\Omega^{p+q}_{\overline{G}}(M)$} and \smash{${\rm d}^{p,q}_0\colon E_0^{p,q}\to E^{p,q+1}_0$} is simply the morphism induced by ${\rm d}$.
\item The $r$-th page is given inductively by
\begin{align*}
E_{r}^{p,q}:={}&\operatorname{Ker}\bigl({\rm d}^{p,q}_{r-1}\bigr)\slash \operatorname{Im}\bigl({\rm d}^{p-r+1,q+r-2}_{r-1}\bigr)\\
={}&\frac{\bigl\{\alpha\in F^{p}_{\mathcal{F}}\Omega^{p+q}_{\overline{G}}(M)\mid {\rm d}\alpha\in F_{\mathcal{F}}^{p+r}\Omega^{p+q+1}_{\overline{G}}(M)\bigr\}}{F^{p+1}_{\mathcal{F}} \Omega^{p+q}_{\overline{G}}(M)+{\rm d}\bigl(F_{\mathcal{F}}^{p-r+1}\Omega_{\overline{G}}^{p+q-1}(M)\bigr)}.
\end{align*}
\item The $r$-th coboundary operator ${\rm d}_r\colon E_r^{p,q}\to E^{p+r,q-r+1}_r$ is again just the map induced by~${\rm d}$ (due to the description of the $r$-th page this has the target specified above and is well defined).
\end{enumerate}
Furthermore, since the filtration is bounded this spectral sequence converges and its final page is isomorphic to the cohomology of the cochain complex \smash{$\Omega^{r}_{\overline{G}}(M)$} known to be isomorphic to the de Rham cohomology of $M$. We call this spectral sequence the spectral sequence of invariant forms. Throughout the rest of the article, the symbol $E_{r}^{p,q}$ will denote this spectral sequence exclusively. The following is a special case of the main result from~\cite{My4} (which we recalled at the beginning of this article as Theorem~\ref{RecG}) which can be in our setting proved verbatim as in~\cite{My3}.
\begin{tw}
Let $(M^{n+s},g)$ be a compact Riemannian manifold with a locally free action of a connected abelian Lie group $G\subset \operatorname{Diff}(M)$ by isometries. Then
\[E_2^{p,q}\cong \bigwedge\nolimits_{H^{p}(M\slash\mathcal{F})}^q\langle\eta_1,\dots,\eta_s\rangle:=H^{p}(M\slash\mathcal{F})\otimes\bigwedge\nolimits^q\langle\eta_1,\dots,\eta_s\rangle.\]
\end{tw}

\begin{proof} Since the operator ${\rm d}$ takes basic forms to basic forms and ${\rm d}\eta_i$ is basic for all $i\in\{1,\dots,s\}$, it is easy to see that ${\rm d}_0$ is in fact equal to the zero operator. Hence, the first page is isomorphic to the $0$-th page.

On the first page by the same observation the operator ${\rm d}_1$ is just the application of ${\rm d}$ to the transverse part of the form (since applying ${\rm d}$ to $\bigwedge ^q\langle\eta_1,\dots,\eta_s\rangle$ decrease $q$ by Proposition~\ref{InvFor}). Hence, the second page is just $H^{p}(M\slash\mathcal{F})\otimes\bigwedge ^q\langle\eta_1,\dots,\eta_s\rangle$.
\end{proof}

\begin{rem} It is apparent that this sequence is a generalization of the sequence presented in~\cite{My3}. Compared to the sequence from~\cite{My4}, we have omitted the adjustment of the metric in Section~3 of the aforementioned paper since it is not necessary in the abelian case. However, it is clear that introducing this modification does not change the resulting sequence (the definition of the sequence itself is independent on the choice of the metric) while the forms $\eta_i$ defined as above and the ones from~\cite{My4} in the case of abelian groups coincide. Moreover, $\bigwedge^q\langle\eta_1,\dots,\eta_s\rangle$ coincides with the Lie algebra cohomology of the Lie algebra $\mathfrak{g}$ of $G$. Consequently, the second page of this sequence is indeed
\[E^{p,q}_2=H^p(M\slash\mathcal{F})\otimes H^{q}(\mathfrak{g}).\]
\end{rem}
We also wish to mention the following consequence of the above discussion which will be used to omit the homological orientability assumption in the final section of the article and is again fully analogous to the corresponding result in~\cite{My3}.

\begin{prop} Let $(M^{n+s},g)$ be a compact oriented Riemannian manifold with a locally free action of a connected abelian Lie group $G\subset \operatorname{Diff}(M)$ by isometries. Then the foliation by orbits is homologically orientable.
\end{prop}

\begin{proof} Without loss of generality, let us assume that $M$ is connected. It is well known (cf.~\cite{E1}) that the top basic cohomology of a~Riemannian foliation on a compact connected manifold is either $0$ or $\mathbb{R}$. In this case, it cannot be~$0$ since then we could compute from the above spectral sequence that \smash{$H^{n+s}_{\rm dR}(M)\cong E^{n,s}_2=0$} which is a contradiction with the orientability of~$M$. Hence, the top basic cohomology is isomorphic to $\mathbb{R}$ which means that the foliation is homologically orientable.
\end{proof}

\section{The degeneration of the spectral sequence}
This section is dedicated to the proof of out main result. The key technical observation is that in the case of an abelian Lie group action the considered spectral sequence becomes somewhat similar to the spectral sequence of a double complex with respect to the computation involved. To be precise, the similarity comes from the fact that with respect to the bigradation the operator~${\rm d}$ has only two potentially non-vanishing parts (the $(1,0)$ and $(2,-1)$ instead of the $(0,1)$ and $(1,0)$ parts). This allows us to track ${\rm d}_r$ using a similar staircase technique as in the case of a~double complex. For convenience we introduce the following notion.

\begin{defi}We will say that the collection $\{{\rm d}\eta_1,\dots,{\rm d}\eta_s\}$ has basic cohomological rank~$k$ if $\dim \langle [{\rm d}\eta_1], \dots, [{\rm d}\eta_s] \rangle = k$ in $H^2(M\slash\mathcal{F})$.
\end{defi}

 Firstly, let us notice that by changing the basis we can choose a basis $\{\tilde{\eta}_1,\dots,\tilde{\eta}_s\}$ for the space $V$ generated (over $\mathbb{R}$) by $\{\eta_1,\dots,\eta_s\}$ such that $\{[{\rm d}\tilde{\eta}_1],\dots,[{\rm d}\tilde{\eta}_k]\}$ are linearly independent and $[{\rm d}\tilde{\eta}_{k+1}]=\cdots =[{\rm d}\tilde{\eta}_s]=0$ (where by $[\alpha]$ we mean the basic cohomology class of~$\alpha$). This can be done by simply considering the map induced by ${\rm d}$ restricted to $V$ with target in $H^2(M\slash\mathcal{F})$ and choosing such a basis $\{\tilde{\eta}_1,\dots,\tilde{\eta}_s\}$ that the last $s-k$ elements span its kernel. Without loss of generality, we will assume that the initial basis $\{\eta_1,\dots,\eta_s\}$ already has this property.

Let us start by treating the case $k=0$. Under this assumption, all of the forms ${\rm d}\eta_i$ are basic exact (i.e., there exist basic forms $\gamma_i$ such that ${\rm d}\gamma_i={\rm d}\eta_i$). This implies immediately that ${\rm d}_2=0$ since ${\rm d}_2(\eta_{i_1}\cdots \eta_{i_l}\alpha)$ is represented by the form
\[\sum (-1)^a\eta_{i_1}\cdots \hat{\eta}_{i_{a}}\cdots \eta_{i_l}{\rm d}\eta_{i_a}\alpha,\]
where the sum runs over $a$. This element is the image through the $(1,0)$ part \big(denoted ${\rm d}^{1,0}$\big) of~${\rm d}$~of
\[\sum \eta_{i_1}\cdots \eta_{i_{a-1}}\gamma_{a}\eta_{i_{a+1}}\cdots \eta_{i_l}\alpha.\]
Note that $\alpha$ was unchanged throughout this process which allows us to perpetuate it for all ${\rm d}_r$. Hence, to compute ${\rm d}_3$ we take the above element (with minus), compute the $(2,-1)$ part of ${\rm d}$ and discern that it is the image through the ${\rm d}^{1,0}$ of
\[-\sum \eta_{i_1}\cdots \eta_{i_{a_1-1}}\gamma_{i_{a_1}}\eta_{i_{a_1+1}}\cdots \eta_{i_{a_2-1}}\gamma_{i_{a_2}}\eta_{i_{a_2+1}}\cdots \eta_{i_l}\alpha,\]
where the sum goes over all $1\leq a_1<a_2\leq l$. Continuing this process inductively will yield that the image of the initial element through ${\rm d}_r$ is the image through ${\rm d}^{1,0}$ of
\[(-1)^r\sum \eta_{i_1}\cdots \eta_{i_{a_1-1}}\gamma_{i_{a_1}}\eta_{i_{a_1+1}}\cdots \eta_{i_{a_r-1}}\gamma_{i_{a_r}}\eta_{i_{a_r+1}}\cdots \eta_{i_l}\alpha,\]
where the sum runs over $1\leq a_1<a_2<\cdots <a_r\leq l$. Hence, indeed the sequence degenerates at the second page in this case.

The plan to generalize this argument for higher $k$ is to keep track of how the basic part changes for ${\rm d}_r$ with $r\leq (k+1)$ and use this along with the above computation to show that given a~representative~$\omega$ of an element of $E_2^{p,q}$ all of the basic forms arising in the chosen representative of ${\rm d}_r(\omega)$ (with $r\geq (k+2)$) are already basic exact under the assumption ${\rm d}_{k+1}\omega=0$. Let us do a practice run by checking the situation for $s=3$ and $k=2$ and a chosen element of the second page represented by $\tilde{\alpha}:=\eta_1\eta_2\eta_3\alpha$ for some basic form $\alpha$. The intention is to showcase the general behaviour of the operator ${\rm d}_k$ before going into detailed computation with an abundance of indices. Firstly, let us note that as before we can assume that the basic class of ${\rm d}\eta_3$ represents zero in basic cohomology and hence there is a basic $1$-form $\gamma_3$ with ${\rm d}\gamma_3={\rm d}\eta_3$. If we now trace the $(2,-1)$ part of ${\rm d}\tilde{\alpha}$, we get the form
\[\eta_2\eta_3{\rm d}\eta_1\alpha-\eta_1\eta_3{\rm d}\eta_2\alpha+\eta_1\eta_2{\rm d}\eta_3\alpha.\]
This is a representative of ${\rm d}_2[\tilde{\alpha}]$. Since in this case we are interested in showing that ${\rm d}_4$ is zero, we assume that the above element represents zero on the second page. Consequently, each of the basic forms ${\rm d}\eta_i\alpha$ are images of the $(1,0)$ part of ${\rm d}$ of some basic form. We denote these forms as $\alpha_i$ and consequently
${\rm d}\eta_i\alpha={\rm d}\alpha_i$.
Moreover, we can put
$\alpha_3=\gamma_3\alpha$.
Hence, we get that
\[{\rm d}^{(1,0)}(\eta_2\eta_3\alpha_1-\eta_1\eta_3\alpha_2+\eta_1\eta_2\gamma_3\alpha)= \eta_2\eta_3{\rm d}\eta_1\alpha-\eta_1\eta_3{\rm d}\eta_2\alpha+\eta_1\eta_2{\rm d}\eta_3\alpha.\]
Continuing ``down the staircase'', we need to compute
\[-{\rm d}^{(2,-1)}(\eta_2\eta_3\alpha_1-\eta_1\eta_3\alpha_2+\eta_1\eta_2\gamma_3\alpha),\]
which gives
\[-[\eta_3({\rm d}\eta_2\alpha_1-{\rm d}\eta_1\alpha_2)+\eta_2({\rm d}\eta_1\gamma_3\alpha-{\rm d}\eta_3\alpha_1)+\eta_1({\rm d}\eta_3\alpha_2-{\rm d}\eta_2\gamma_3\alpha)].\]
Trying to find the preimage of this element through the $(1,0)$ part, we get that
\begin{gather*}
{\rm d}(-\gamma_3\alpha_1)={\rm d}\eta_1\gamma_3\alpha-{\rm d}\eta_3\alpha_1,\qquad
{\rm d}(\gamma_3\alpha_2)={\rm d}\eta_3\alpha_2-{\rm d}\eta_2\gamma_3\alpha.
\end{gather*}
Hence, by assuming as before that ${\rm d}_3[\tilde{\alpha}]=0$ we get one new assumption. Namely, that $\eta_3({\rm d}\eta_2\alpha_1-{\rm d}\eta_1\alpha_2)$ is a combination of an element in the image of the $(2,-1)$ and $(1,0)$ parts of~${\rm d}$ of two different forms (of appropriate bi-degree). Hence, by changing $\alpha_1$ and $\alpha_2$ if need be we can assume that $({\rm d}\eta_2\alpha_1-{\rm d}\eta_1\alpha_2)$ is an image through the $(1,0)$ part of ${\rm d}$ of some basic form~$\alpha_{12}$. Hence, continuing the process using these assumptions we get that ${\rm d}_4[\tilde{\alpha}]$ is represented by
\[-[{\rm d}\eta_3\alpha_{12}-{\rm d}\eta_2\gamma_3\alpha_1+{\rm d}\eta_1\gamma_3\alpha_2],\]
which is a basic exact form since
\[{\rm d}(\gamma_3\alpha_{12})={\rm d}\eta_3\alpha_{12}-{\rm d}\eta_2\gamma_3\alpha_1+{\rm d}\eta_1\gamma_3\alpha_2.\]
Hence , ${\rm d}_4[\tilde{\alpha}]=0$ and necessarily by degree reasons all subsequent ${\rm d}_r[\tilde{\alpha}]=0$. From the above example, we make the following observations:
\begin{itemize}\itemsep=0pt
\item The preimage of elements with ${\rm d}\eta_i$ can be (in the proof above and as we will later see in the following computation) found based on the assumption that the given representative induces an element on the $(k+1)$-th page.
\item The preimages of elements with ${\rm d}\eta_i$ which is basic exact can be taken without changing the rest of the basic form.
\item The reasoning presented above gives us a pretty good guess on how the representatives of each ${\rm d}_r[\tilde{\alpha}]$ should look.
\end{itemize}
This leads us to the following notation which will be useful in the proper computation:
\begin{Not} Let $\alpha$ denote a homogenous differential form $($with respect to the bigrading$)$ which represents a given class on the second page. We can split this form in a unique way as a~sum
\[\alpha=\sum \eta_{i_1}\cdots \eta_{i_a}\Bigl(\sum \eta_{j_1}\cdots \eta_{j_b}\alpha_{i_1,\dots,i_a}^{j_1,\dots,j_b}\Bigr),\]
where
\begin{enumerate}\itemsep=0pt
\item[$1.$] The first sum goes over all $0 \leq a\leq s-k$ and $k<i_1<\cdots <i_a\leq s$.
\item[$2.$] In the second sum, $b$ is fixed so that $a+b$ give the second part of the bidegree and the sum goes over all $0<j_1<\cdots <j_b\leq k$.
\item[$3.$] The forms \smash{$\alpha_{i_1,\dots,i_a}^{j_1,\dots,j_b}$} are basic and closed since $\alpha$ induces an element of the second page.
\end{enumerate}
Recall that we do this under the assumption that $\eta_{i}$ with $k<i\leq s$ are basic exact and $\{{\rm d}\eta_1,\dots,{\rm d}\eta_s\}$ have basic cohomological rank~$k$. Hence, for each $\eta_i$ with $k<i\leq s$ we denote by~$\gamma_i$ the basic $1$-form with ${\rm d}\eta_i={\rm d}\gamma_i$. The reason for splitting the indices like this is in order to track how the basic forms change as we go down the staircase. Hence, we define~\smash{$\alpha_{i_1,\dots,i_a}^{j_1,\dots,j_{b-1}}$} as follows. Firstly, compute ${\rm d}^{(2,-1)}\alpha$ and consider the basic form $\beta$ corresponding to $\eta_{i_1}\cdots \eta_{i_a}\eta_{j_1}\cdots \eta_{j_{b-1}}$. When choosing an element in the preimage for this form through ${\rm d}^{(1,0)}$ take the components coming from changing exact $\eta_i$ to ${\rm d}\eta_i$ and find its preimage through ${\rm d}^{(1,0)}$ as described in the exact case $($this simply changes this copy of ${\rm d}\eta_i$ to $\gamma_i$ at this step$)$. Then take an element of the preimage of the remainder of $\beta$ arbitrarily $($this will be later adjusted in the further steps of the construction$)$. This element is \smash{$\alpha_{i_1,\dots,i_a}^{j_1,\dots,j_{b-1}}$}.

 We continue this construction inductively defining all \smash{$\alpha_{i_1,\dots,i_a}^{j_1,\dots,j_{b-l}}$} with an arbitrary number of upper indexes missing by repeating it for resulting element from the previous step $($note that we do not change the signs here$)$. This can be done since we assume that ${\rm d}_r[\alpha]=0$ for $2\leq r\leq (k+1)$ and hence all the necessary preimages are indeed non-empty.
\end{Not}

Let us now use this notation to show that ${\rm d}_r[\alpha]=0$ for $r\geq (k+2)$. We do this by simply listing all the chosen elements of the preimages through ${\rm d}^{(1,0)}$ which occur as we go down the staircase. Hence, we get the first element
\[\sum \eta_{i_1}\cdots \gamma_{i_c}\cdots \eta_{i_a}\Bigl(\sum \eta_{j_1}\cdots \eta_{j_b}\alpha_{i_1,\dots,i_a}^{j_1,\dots,j_b}\Bigr)+\sum\eta_{i_1}\cdots \eta_{i_a}\Bigl(\sum \eta_{j_1}\cdots \eta_{j_{b-1}}\alpha_{i_1,\dots,i_a}^{j_1,\dots,j_{b-1}}\Bigr),\]
where the first sum goes additionally over all possible values of~$c$. Similarly as the above element has two sums the element of the preimage through ${\rm d}^{1,0}$ of the form representing ${\rm d}_3[\alpha]$ will have three elements (depending
on how many upper indices are omitted in \smash{$\alpha_{i_1,\dots,i_a}^{j_1,\dots,j_b}$}). These are
\begin{gather*}
-\sum \eta_{i_1}\cdots \gamma_{i_{c_1}}\cdots \gamma_{i_{c_2}}\cdots \eta_{i_a}\Bigl(\sum \eta_{j_1}\cdots \eta_{j_b}\alpha_{i_1,\dots,i_a}^{j_1,\dots,j_b}\Bigr),\\
-\sum \eta_{i_1}\cdots \gamma_{i_{c_1}}\cdots \eta_{i_a}\Bigl(\sum \eta_{j_1}\cdots \eta_{j_{b-1}}\alpha_{i_1,\dots,i_a}^{j_1,\dots,j_{b-1}}\Bigr),\\
-\sum \eta_{i_1}\cdots \eta_{i_a}\Bigl(\sum \eta_{j_1}\cdots \eta_{j_{b-2}}\alpha_{i_1,\dots,i_a}^{j_1,\dots,j_{b-2}}\Bigr),
\end{gather*}
where in the first expression the sum goes additionally over $c_1<c_2$. Note that some of the terms above might be null due to lack of appropriate indices (e.g., in the example above the top term of these three was already zero). Iterating this process even further (and remembering to change the sign), we get the preimage of the form representing ${\rm d}_4[\alpha]$ as the sum of the expressions
\begin{gather*}
\sum \eta_{i_1}\cdots \gamma_{i_{c_1}}\cdots \gamma_{i_{c_2}}\cdots \gamma_{i_{c_3}}\cdots \eta_{i_a}\Bigl(\sum \eta_{j_1}\cdots \eta_{j_b}\alpha_{i_1,\dots,i_a}^{j_1,\dots,j_b}\Bigr),\\
\sum \eta_{i_1}\cdots \gamma_{i_{c_1}}\cdots \gamma_{i_{c_2}}\cdots \eta_{i_a}\Bigl(\sum \eta_{j_1}\cdots \eta_{j_{b-1}}\alpha_{i_1,\dots,i_a}^{j_1,\dots,j_{b-1}}\Bigr),\\
\sum \eta_{i_1}\cdots \gamma_{i_{c_1}}\cdots \eta_{i_a}\Bigl(\sum \eta_{j_1}\cdots \eta_{j_{b-2}}\alpha_{i_1,\dots,i_a}^{j_1,\dots,j_{b-2}}\Bigr),\\
\sum \eta_{i_1}\cdots \eta_{i_a}\Bigl(\sum \eta_{j_1}\cdots \eta_{j_{b-3}}\alpha_{i_1,\dots,i_a}^{j_1,\dots,j_{b-3}}\Bigr).
\end{gather*}
Note that as before some of these elements might be zero (e.g., in the example above all the elements except for the third where already zero). The computation continues as above through all ${\rm d}_r[\alpha]$ until all the $\eta_i$ are eliminated. Let us list the final two elements of this. Then the penultimate element consists of two sums
\[\pm\sum \gamma_{i_1}\cdots \eta_{i_{c}}\cdots \gamma_{i_{a}}(\alpha_{i_1,\dots,i_{a}}),\qquad \pm\sum \gamma_{i_1}\cdots \gamma_{i_{a}}\Bigl(\sum \eta_{j}\alpha_{i_1,\dots,i_{a}}^{j}\Bigr).\]
The final element of this process is just
\[\mp\sum \gamma_{i_1}\cdots \gamma_{i_{a}}\alpha_{i_1,\dots,i_{a}}.\]
Hence, we have shown that if ${\rm d}_r[\alpha]=0$ for $r\leq k+1$ then it has to be zero for all $r$ which finishes the proof of our main result.

\section{Examples}
In this section, we give a number of examples in order to show that the estimates from the previous section are sharp. More than that we show that for any $s\in\mathbb{N}$, $0\leq k\leq s$ and~$2\leq l\leq k+2$ there is a compact Riemannian manifold $(M^{n+s},g)$ with an action of a connected $s$-dimensional abelian Lie group $G\subset \operatorname{Diff}(M)$ by isometries such that
\begin{enumerate}\itemsep=0pt
\item The set of elements $\{{\rm d}\eta_1,\dots,{\rm d}\eta_s\}$ has basic cohomological rank $k$.
\item The spectral sequence degenerates on the $l$-th page.
\end{enumerate}
We start with a $2s$-torus $\mathbb{T}^{2s}$ and name the standard representatives of generators of its first cohomology as $\{\overline{\alpha}_1,\dots,\overline{\alpha}_s,\overline{\beta}_1,\dots,\overline{\beta}_{s-1},\overline{\sigma}\}$. Next we are going to consider principal circle bundles over this torus with appropriate curvatures by using the following result which we quote verbatim from~\cite{BG}.

\begin{tw}[{\cite[Theorem 7.1.6]{BG}}] Let $(Z; \omega; J)$ be an almost K\"ahler orbifold with $[p^*\omega]\in H^2_{\rm orb}(Z,\mathbb{Z})$ and let $M$ denote the total space of the circle $V$-bundle defined by the class $[\omega]$. Then the orbifold $M$ admits a $K$-contact structure $(\xi; \Phi;\eta; g)$ such that ${\rm d}\eta = \pi^*\omega$, where $\pi\colon M\to Z$ is the natural orbifold projection map. Furthermore, if all the local uniformizing groups of $Z$ inject into the structure group $S^1$, then M is a smooth $K$-contact manifold.
\end{tw}

As we do not want to go too deep into the language used in this theorem, we explicitly state an immediate corollary of this theorem which is relevant to our construction below. This is just a restriction of the statement to the case of smooth manifolds (omitting the language of orbifolds) and separating the part of the thesis relevant to us.

\begin{cor} For an almost K\"ahler manifold $(N,\omega,J)$ with a principal circle bundle $($with total space $M)$ represented by the class $[\omega]\in H^2_{\rm dR}(M)$ which is integral there exists a Riemannian metric $g$ on $M$, satisfying ${\rm d}\eta=\pi^*\omega$ $($where $\eta(\bullet)=g(\xi,\bullet)$ and $\xi$ is a chosen vector field representing the $\mathbb{S}^1$-action$)$.
\end{cor}
We use this as follows. Consider the projection $\pi_{\alpha_i}$ (resp.\ $\pi_{\beta_i}$) of the given torus $\mathbb{T}^{2s}$ onto the torus $\mathbb{T}^2$ corresponding to $\{\overline{\alpha}_i,\overline{\sigma}\}$ \smash{\big(resp.\ $\bigl\{\overline{\beta}_i,\overline{\sigma}\bigr\}$\big)}. We will denote the forms corresponding to~$\{\overline{\alpha}_i,\overline{\sigma}\}$ \smash{\big(resp.\ $\bigl\{\overline{\beta}_i,\overline{\sigma}\bigr\}$\big)} on this torus by $\{\tilde{\alpha}_i,\sigma_{\alpha_i}\}$ \smash{\big(resp.\ $\bigl\{\tilde{\beta}_i,\sigma_{\beta_i}\bigr\}$\big)}. Since such a torus is K\"ahler with symplectic form $\omega_{\alpha_i}=\sigma_{\alpha_i}\wedge\tilde{\alpha}_i$ \smash{\big(resp.\ $\omega_{\beta_i}=\tilde{\beta}_i\wedge\sigma_{\beta_i}$\big)} we can use the above corollary to consider principal bundles $p_{\alpha_i}\colon \tilde{P}_{\alpha_i}\to\mathbb{T}^2$ and $p_{\beta_i}\colon \tilde{P}_{\beta_i}\to\mathbb{T}^2$ with corresponding one forms~$\tilde{\eta}_i$ and~$\tilde{\gamma}_i$ with \smash{${\rm d}\tilde{\eta}_i=p^*_{\alpha_i}\omega_{\alpha_i}$} and \smash{${\rm d}\tilde{\gamma}_i=p_{\beta_i}^*\omega_{\beta_i}$}. Finally, we can pull these bundles back by the respective projections to $\mathbb{T}^{2s}$ to get bundles \smash{$P_{\alpha_i}=\pi^*_{\alpha_i}\tilde{P}_{\alpha_i}$} and \smash{$P_{\beta_i}=\pi^*_{\beta_i}\tilde{P}_{\beta_i}$} and consider ${\pi\colon M\to\mathbb{T}^{2s}}$ to be the direct sum of all these circle bundles. Due to the construction, we additionally have the obvious maps $\tilde{\pi}_{\alpha_i}\colon P_{\alpha_i}\to\tilde{P}_{\alpha_i}$, $\tilde{\pi}_{\beta_i}\colon P_{\beta_i}\to\tilde{P}_{\beta_i}$, $f_{\alpha_i}\colon M\to P_{\alpha_i}$ and $f_{\beta_i}\colon M\to P_{\beta_i}$.

Under the above construction, the cotangent bundle of $M$ is trivialized by the sections \smash{$\alpha_i:=\pi^*\overline{\alpha}_i$}, \smash{$\beta_i:=\pi^*\overline{\beta}_i$}, \smash{$\sigma:=\pi^*\overline{\sigma}$}, \smash{$\eta_i:=f_{\alpha_i}^*\circ\tilde{\pi}^*_{\alpha_i}\tilde{\eta_i}$} and \smash{$\gamma_i:=f_{\beta_i}^*\circ\tilde{\pi}^*_{\beta_i}\tilde{\gamma_i}$}. We define a Riemannian metric on $M$ by demanding that this trivialization is orthonormal. Moreover, since $M$ is a~$\mathbb{T}^{2s-1}$ principal bundle over $\mathbb{T}^{2s}$, it admits an action of $\mathbb{T}^{2s-1}$. Let $\{\xi_1,\dots,\xi_s\}$ be the vector fields corresponding to the bundles $P_{\alpha_i}$ in the direct sums and hence inducing an action of~${\mathbb{T}^s\subset\mathbb{T}^{2s-1}}$ on $M$. These vector fields are dual with respect to the above metric to the forms~$\eta_i$. Finally, let us describe the operator ${\rm d}$ with respect to this frame
${\rm d}\alpha_i={\rm d}\beta_i={\rm d}\sigma=0$, ${\rm d}\eta_i=\sigma\wedge\alpha_i$, ${\rm d}\gamma_i=\beta_i\wedge\sigma$.
We will show that for this example (taken with the $\mathbb{T}^s$-action described above) the spectral sequence degenerates precisely at \smash{$E^{p,q}_{s+2}$}. To show this we will prove that the element represented by $\eta_1\cdots \eta_s\beta_1\cdots \beta_{s-1}$ on the second page survives to $E^{p,q}_{s+1}$ and that ${\rm d}_{s+1}(\eta_1\cdots \eta_s\beta_1\cdots \beta_{s-1})\neq 0$.

Firstly, we need to compute the basic cohomology of the foliation~$\mathcal{F}$ by orbits. For this, notice that the $\mathbb{T}^s$-action forms a $\mathbb{T}^s$ principal bundle which can be described as the pullback of the direct sum of $P_{\alpha_i}$ to the total space $\tilde{M}$ of the direct sum of $P_{\beta_i}$. Consequently, the basic cohomology is just the cohomology of~$\tilde{M}$. The cohomology of~$\tilde{M}$ can be easily computed from the above description and the considered spectral sequence. For this, one notes that the second page is just the (appropriately graded) real exterior algebra on the vector space generated by $\{\alpha_1,\dots,\alpha_s,\beta_1,\dots,\beta_{s-1},\sigma,\gamma_1,\dots,\gamma_{s-1}\}$, where the basic part is represented by the first $2s$-elements. It is also immediate that the kernel of ${\rm d}_2$ consists of sums of elements which are either multiplicities of $\sigma$ or are of the form
\[(\gamma_{i_1}\beta_{i_1})\cdots (\gamma_{i_k}\beta_{i_k})(\gamma_{j_1}\beta_{h_1}+\gamma_{h_1}\beta_{j_1})\cdots (\gamma_{j_l}\beta_{h_l}+\gamma_{h_l}\beta_{j_l})\omega,\]
where $\omega$ is a product of $\alpha_i$ and $\beta_i$. On the other hand, the image is just described by computing the elements ${\rm d}(\gamma_{i_1}\cdots \gamma_{i_k})$ which amounts to identifying $-\sigma\gamma_{i_1}(\sum \gamma_{i_2}\cdots \gamma_{i_{j-1}}\beta_{i_j}\gamma_{i_{j+1}}\cdots \gamma_{i_k})\omega$ with $\sigma\beta_{i_1}\gamma_{i_2}\cdots \gamma_{i_k}\omega$ (where~$\omega$ is as before). Moreover, the consequent description of the third page can be taken as the description of the cohomology of~$\tilde{M}$ since this spectral sequence degenerates at the third page. To see this, it suffices to note that for the chosen representatives of the classes on the second page an element is in the kernel of ${\rm d}_2$ if and only if the corresponding form is closed. Consequently, for any element of degree $k$ which lives to the third page, we can extend it by $0$ in other bidegrees $(p,q)$ with $p+q=k$ to get the corresponding representative of the cohomology class in~$\tilde{M}$.

Let us now show that the class of the element $\eta_1\cdots \eta_s\beta_1\cdots \beta_{s-1}$ is in the kernels of ${\rm d}_{k}$ for~${1\leq k\leq s}$. To compute this, let us note again that in the case of an abelian Lie group action the considered spectral sequence becomes somewhat similar to the spectral sequence of a double complex with respect to the computation involved. This allows us to track ${\rm d}_k(\eta_1\cdots \eta_s\beta_1\cdots \beta_{s-1})$ by using the same staircase technique as in the case of a double complex. Hence, we start by noting that this indeed gives an element of the second page since ${\rm d}(\beta_1\cdots \beta_{s-1})=0$. Next, one can easily compute that ${\rm d}_2(\eta_1\cdots \eta_s\beta_1\cdots \beta_{s-1})$ is represented by the class of
\[\sum (-1)^i\eta_1\cdots \eta_{i-1}\sigma\alpha_i\eta_{i+1}\cdots \eta_s\beta_1\cdots \beta_{s-1},\]
we can now alter each term in the above sum by pushing $\sigma\alpha_i$ past $\beta_1$ and then pushing the newly constructed $\beta_1\sigma$ to the front of the term (since all permutations are done by changing the place of the $2$-form this does not change the sign). Hence, we arrive at
\[\sum (-1)^i\beta_1\sigma\eta_1\cdots \eta_{i-1}\eta_{i+1}\cdots \eta_s\alpha_i\beta_2\cdots \beta_{s-1},\]
which is the image of
\[\sum (-1)^i\gamma_1\eta_1\cdots \eta_{i-1}\eta_{i+1}\cdots \eta_s\alpha_i\beta_2\cdots \beta_{s-1},\]
through the $(1,0)$ part of ${\rm d}$. Hence, indeed ${\rm d}_2(\eta_1\cdots \eta_s\beta_1\cdots \beta_{s-1})$ is zero in $E^{p,q}_2$. Unfortunately, to simplify further computation we need to modify this form slightly within the preimage of the given element through ${\rm d}^{(1,0)}$. In the above expression, we have elected to use $\sigma$ in every term to transform $\beta_1$ to $\gamma_1$. However, a similar computation can be conducted for any $\beta_j$. Averaging the resulting forms for each $\beta_j$ will give
\[\frac{1}{s-1}\sum (-1)^{i+j-1}\gamma_j\eta_1\cdots \eta_{i-1}\eta_{i+1}\cdots \eta_s\alpha_i\beta_1\cdots \beta_{j-1}\beta_{j+1}\cdots \beta_{s-1},\]
We can then continue the computation by computing $(2,-1)$ part of ${\rm d}$ applied to the above form (with a minus sign) and identifying a representative of its preimage through the $(1,0)$ part. Hence, similarly as above, we get{\samepage
\begin{gather*}
\begin{split}
& -\frac{1}{{s-1\choose 2}}\sum (-1)^{i_1+j_1+i_2+j_2-2}\gamma_{j_1}\gamma_{j_2}\eta_1\cdots \eta_{i_1-1}\eta_{i_1+1}\cdots \eta_{i_2-1}\eta_{i_2+1}\cdots\eta_s\\
& \hphantom{-\frac{1}{{s-1\choose 2}}\sum }{}\wedge
\alpha_{i_1}\alpha_{i_2}\beta_1\cdots \beta_{j_1-1}\beta_{j_1+1}\cdots \beta_{j_2-1}\beta_{j_2+1}\cdots \beta_{s-1},
\end{split}
\end{gather*}
where the sum goes over all $i_1<i_2$ and $j_1<j_2$.}

\begin{rem} It is perhaps somewhat instructive to see that the signs indeed agree between the two terms adding up to
\begin{gather*}
-\frac{2}{(s-1)(s-2)} (-1)^{i_1+j_1+i_2+j_2-2}\gamma_{j_1}\gamma_{j_2}\eta_1\cdots \eta_{i_1-1}\eta_{i_1+1}\cdots \eta_{i_2-1}\eta_{i_2+1}\cdots\eta_s \\
\qquad{}\wedge
\alpha_{i_1}\alpha_{i_2}\beta_1\cdots \beta_{j_1-1}\beta_{j_1+1}\cdots \beta_{j_2-1}\beta_{j_2+1}\cdots \beta_{s-1}.
\end{gather*}
Firstly, let us note that if $i_1<i_2$ then the computation of the sign is exactly the same as above. On the other hand, if $i_2<i_1$, then we need to change the sign twice. Firstly, since now the parity of the position of $\eta_{i_2}$ has changed and secondly to invert $\alpha_{i_1}$ with~$\alpha_{i_2}$.
\end{rem}
The process now continues (with the coefficients in the $k$-th step up to sign equal to \smash{$1/{s-1\choose k}$}) to finally arrive at the fact that ${\rm d}_s(\eta_1\cdots \eta_s\beta_1\cdots \beta_{s-1})$ indeed represents the zero class since it is the image through ${\rm d}^{(1,0)}$ of the form
\[\pm\sum (-1)^{\frac{1}{2}(s-1)s-i+\frac{1}{2}(s-2)(s-1)-s}\gamma_1\cdots \gamma_{s-1}\eta_i\alpha_{1}\cdots \alpha_{i-1}\alpha_{i+1}\cdots \alpha_s.\]

This not only finishes the proof that $\eta_1\cdots \eta_s\beta_1\cdots \beta_{s-1}$ indeed represents an element of $E_{s+1}^{p,q}$ but also is a starting point for proving that ${\rm d}_{s+1}(\eta_1\cdots \eta_s\beta_1\cdots \beta_{s-1})\neq 0$ since it allows us to identify the representative of this class as
\[\mp (-1)^{\frac{1}{2}(s-1)s+\frac{1}{2}(s-2)(s-1)}s\gamma_1\cdots \gamma_{s-1}\sigma\alpha_1\cdots \alpha_s.\]
To see that this represents a nonzero class on the final page, let us note that from the description of the cohomology of~$\tilde{M}$ any class in it is represented by an $\mathbb{R}$-linear combination of elements from the chosen basis of the exterior product. Consequently, we can do the same for any element of the second page of the spectral sequence of the $\mathbb{T}^s$ action on~$M$. Hence, it suffices to show that written in such a basis no image of an element contains the above term. To see this, let us first notice that terms containing $\gamma_i$ are not in the image through ${\rm d}$ on $\tilde{M}$ of any element that does not contain~$\gamma_i$ already. Hence, $\gamma_i$ can be only produced one at a time as we apply the $(2,-1)$ part of~${\rm d}$. Hence, if this class was a term in an image through ${\rm d}_k$ (for $k\leq s$) of some element then that element written in this basis would have to have a term already containing at least one~$\gamma_i$ (in fact, it would have to have $s-k+1$ such forms in any term that contributes to the final term). By the description of representatives of classes in the cohomology of $\tilde{M}$, it can be paired either into components with $\sigma$, $\gamma_i\beta_i$ or $\gamma_i\beta_j+\gamma_j\beta_i$. However, the component with $\sigma$ is automatically closed and hence its image via ${\rm d}^{2,-1}$ will be zero. Moreover, due to the limited options, we are able quite accurately to predict how the staircase argument will change the given element without~$\sigma$. Hence, since our end goal has no $\beta_i$ each such component would have to change to~$\gamma_i$ at some point in at least one term. This allows us to conclude that elements with $\beta_i\gamma_i$ cannot contribute to the end goal (since at some point the staircase argument for such an element would give a multiplicity of $\gamma_i\gamma_i$ after pulling back through ${\rm d}^{1,0}$ which gives a contradiction). Similarly, one can see that the terms
$\gamma_i\beta_j+\gamma_j\beta_i$ cannot themselves be pulled back (since they would give $\gamma_i\gamma_j+\gamma_j\gamma_i$) which gives grounds for eliminating the final type of elements. More precisely, given an element with component $\gamma_i\beta_j+\gamma_j\beta_i$ we can again effortlessly see to what element the elements arising from splitting this $2$-form will contribute in the end. Such elements will have all $\beta_i$ replaced by $\gamma_i$, and all $\eta_i$ replaced by $\alpha_i$. Taking into account the arising coefficient (which can be computed similarly as in the previous paragraph) and taking an average over all the components that contribute to the change of the given component, we conclude that if this process can be conducted then it simply changes $\beta_i$ to $\gamma_i$ up to a set constant. Moreover, the terms with $\gamma_i\beta_j$ and $\gamma_j\beta_i$ will undergo the same process (after we gather all the necessary forms for pulling them back at each step) in the end and hence the coefficients standing next to them will differ by a minus sign. Consequently, a form consisting of elements of this form if it could be taken down the staircase to the final row would contribute zero to the term in question.

Hence, indeed the form
$\mp (-1)^{\frac{1}{2}(s-1)s+\frac{1}{2}(s-2)(s-1)}s\gamma_1\cdots \gamma_{s-1}\sigma\alpha_1\cdots \alpha_s$
is not a component in any element of the image of ${\rm d}_k$ and consequently represents a nonvanishing element in $E_{s+1}^{p,q}$ proving that this sequence degenerates at the $s+2$ page.

Let us finish this section by modifying these examples slightly so that the page on which the sequence degenerates is independent (except for the obvious limitations and the main result of the previous section) on the dimension of the group or on the basic cohomological rank of $\{{\rm d}\eta_1,\dots,{\rm d}\eta_s\}$. This is done by simply taking products with appropriate spaces. To increase the dimension of the group without increasing the basic cohomological rank of $\{{\rm d}\eta_1,\dots,{\rm d}\eta_s\}$ it suffices to take a product with a circle acting on itself. By a similar trick, we can increase the dimension of the group simultaneously with the basic cohomological rank of $\{{\rm d}\eta_1,\dots,{\rm d}\eta_s\}$ by taking a product with a principal bundle isomorphic to one of the bundles $\tilde{P}_{\alpha_i}$. It is easy to see that making a product with any number of the above spaces will not change the above computation in a significant way allowing us to conclude that the sequence of a product of one of the above examples with a number of the above spaces will still degenerate at the same page.

\section{Relation to harmonic forms}
In this section, we study some applications of the above sequence. In doing so, we underline the usefulness of an upper bound on the number of page on which the sequence degenerates. As already mentioned in the introduction, some of the motivation for this work came from the study of $\mathcal{K}$-structures (cf.~\cite{B}) and $q$-contact manifolds (cf.~\cite{F}). In particular, it is evident that we can use our spectral sequence to study cohomological properties of such manifolds. However, due to a more geometric nature of such manifolds, we wish to give a more geometrically significant application as well. Hence, we present a possible additional use for this sequence in calculating the harmonic forms of a given manifold in terms of basic harmonic forms. Unfortunately, due to the vast variety of possible interplay between the forms ${\rm d}\eta_i$, it seems that a clear cut classification as was done for $\mathcal{S}$-manifolds in~\cite{My4} is nearly impossible in full generality. However, we can still apply the sequence in a similar fashion to compute the harmonic forms in a given example.

There are two key observations which make this approach viable. The first one is that we do not exclude harmonic forms by passing to invariant forms. This is covered by the following well-known result.

\begin{prop} If $\alpha$ is a harmonic form and $X$ is a Killing vector field, then $\mathcal{L}_X\alpha=0$. In particular, if a connected Lie group $G$ acts on $M$ by isometries, then harmonic forms are invariant with respect to this action.
\end{prop}

The second key observation allowing us to apply this line of thought effectively is that as the spectral sequence allows us to compute the cohomology of~$M$ using the basic cohomology of the foliation by orbits and the Lie algebra cohomology of the acting abelian group the Hodge star operator respects this separation. More precisely, by straightforward computation, we get the following.

\begin{prop}
Let $(M^{n+s},g)$ be a compact oriented Riemannian manifold with a locally free isometric action of a connected abelian Lie group $G$ $($of dimension $s)$ and $\{\eta_1,\dots,\eta_s\}$ are pointwise orthonormal. Moreover, let $i=(i_1,\dots,i_k)$ be an ordered subset of $\{1,\dots,s\}$ with complement $j=(j_1,\dots,j_{s-k})$. Then the following relation between the Hodge star operator $*$ and the basic Hodge star operator $*_{b}$ holds for any transverse $r$-form $\alpha$:
\[*(\eta_{i_1}\cdots \eta_{i_k}\alpha)={\rm sign}(i_1,\dots,i_k,j_1,\dots,j_{s-k})(-1)^{(s-k)r}\eta_{j_1}\cdots \eta_{j_{s-k}}*_b\alpha.\]
\end{prop}
\begin{proof} It suffices to note that under the above assumptions $\eta_i$ are pointwise orthogonal to any $1$-form which vanishes on $T\mathcal{F}$ and then the proposition follows by straightforward pointwise computation.
\end{proof}

Let us see how we can employ the above observations to finding the harmonic forms in one of the examples from the previous section.
\begin{egg} Here we will compute explicitly the cohomology and basic harmonic forms (with respect to the appropriate metric) of one of the examples from the previous section. Let us pick the example with $k=s=2$. Consequently, we have a $7$-manifold $M$ with trivial cotangent bundle spanned by the forms $\{\alpha_1,\alpha_2,\beta,\sigma,\gamma,\eta_1,\eta_2\}$ (which we use to define the metric by demanding that this collection is orthonormal) which are subjugate to the relations
${\rm d}\alpha_1={\rm d}\alpha_2={\rm d}\sigma={\rm d}\beta=0$, ${\rm d}\gamma=\beta\sigma$,
${\rm d}\eta_1=\sigma\alpha_1$, ${\rm d}\eta_2=\sigma\alpha_2$.
Moreover, we recall that $M$ is a torus bundle over a $5$-manifold $\tilde{M}$ which is a circle bundle. Let us start by computing the cohomology of~$\tilde{M}$ using our spectral sequence. With the given presentation of the second page (which amounts to it being an exterior algebra on $\langle\alpha_1,\alpha_2,\beta,\sigma,\gamma\rangle$) and the fact that this sequence has only two rows, it suffices to compute ${\rm d}^{(2,-1)}$ which simply changes~$\gamma$ (or more precisely a form on~\smash{$\tilde{M}$} which is pulled back to $\gamma$ on $M$) to ${\rm d}\gamma=\beta\sigma$. Hence, the second page of this sequence is of the form
\[\begin{tabular}{|c c c c c}
$\mathbb{R}$ & $\mathbb{R}^4$ & $\mathbb{R}^6$&$\mathbb{R}^4$& $\mathbb{R}$\\
$\mathbb{R}$ & $\mathbb{R}^4$ & $\mathbb{R}^6$&$\mathbb{R}^4$& $\mathbb{R}$
\\\hline
\end{tabular}\]
whereas the third page after the necessary computation is of the form
\[\begin{tabular}{|c c c c c}
$0$ & $\mathbb{R}^2$ & $\mathbb{R}^5$&$\mathbb{R}^4$& $\mathbb{R}$\\
$\mathbb{R}$ & $\mathbb{R}^4$ & $\mathbb{R}^5$&$\mathbb{R}^2$& $0$
\\\hline
\end{tabular}\]
The generators in terms of $\{\alpha_1,\alpha_2,\beta,\sigma,\gamma\}$ are simple enough to compute and we will invoke them in further computations when needed. By taking direct sums on the diagonal, we get the cohomology of $\tilde{M}$. We now use this to compute the cohomology of $M$. The second page of the spectral sequence related to the $\mathbb{T}^2$ action on~$M$ is thus
\[\begin{tabular}{|c c c c c c}
$\mathbb{R}$ & $\mathbb{R}^4$ & $\mathbb{R}^7$&$\mathbb{R}^7$& $\mathbb{R}^4$&$\mathbb{R}$\\
$\mathbb{R}^2$ & $\mathbb{R}^{8}$ & $\mathbb{R}^{14}$&$\mathbb{R}^{14}$& $\mathbb{R}^{8}$&$\mathbb{R}^2$\\
$\mathbb{R}$ & $\mathbb{R}^4$ & $\mathbb{R}^7$&$\mathbb{R}^7$& $\mathbb{R}^4$&$\mathbb{R}$
\\\hline
\end{tabular}\]
A somewhat more difficult computation gives us the third page. For the sake of brevity, we will not show it here in full (since computing any one position on the page boils down to elementary computation). Moreover, some further light on these computations should be shed by the later description of harmonic forms. The third page is then as follows:
\[\begin{tabular}{|c c c c c c}
$0$ & $\mathbb{R}^2$ & $\mathbb{R}^6$&$\mathbb{R}^5$& $\mathbb{R}^4$&$\mathbb{R}$\\
$0$ & $\mathbb{R}^{7}$ & $\mathbb{R}^{11}$&$\mathbb{R}^{11}$& $\mathbb{R}^{7}$&$0$\\
$\mathbb{R}$ & $\mathbb{R}^4$ & $\mathbb{R}^5$&$\mathbb{R}^6$& $\mathbb{R}^2$&$0$
\\\hline
\end{tabular}\]
Finally, as mentioned in the description of these examples in the previous section, ${\rm d}_3$ does not vanish on the form $\eta_1\eta_2\beta$. By computing that it indeed vanishes on the other generator (namely~$\eta_1\eta_2\sigma$) in that position and by noting that due to degree reasons this is the only position on which ${\rm d}_3$ might not be zero, we get the final page
\[\begin{tabular}{|c c c c c c}
$0$ & $\mathbb{R}$ & $\mathbb{R}^6$&$\mathbb{R}^5$& $\mathbb{R}^4$&$\mathbb{R}$\\
$0$ & $\mathbb{R}^{7}$ & $\mathbb{R}^{11}$&$\mathbb{R}^{11}$& $\mathbb{R}^{7}$&$0$\\
$\mathbb{R}$ & $\mathbb{R}^4$ & $\mathbb{R}^5$&$\mathbb{R}^6$& $\mathbb{R}$&$0$
\\\hline
\end{tabular}\]
Hence, the computation of the cohomology of $M$ is indeed concluded. To find the harmonic forms, it suffices to find the representatives of each of the classes keeping in mind that the behaviour of the Hodge star respects the splitting we have had on the second page. Let us write them down and then comment on the computation (verification that these forms are indeed harmonic is straightforward). For example, the harmonic forms corresponding to the lower row in each degree are spanned by the following linearly independent sets:
\begin{gather*}
\{1\},\qquad
\{\alpha_1,\alpha_2,\beta,\sigma\},\qquad
\{\alpha_1\beta,\alpha_2\beta,\alpha_1\alpha_2,\gamma\sigma,\gamma\beta\},\\
\{\alpha_2\sigma\gamma-\eta_2\beta\sigma,\alpha_1\sigma\gamma-\eta_1\beta\sigma, \beta\sigma\gamma,\alpha_1\alpha_2\beta,\alpha_2\beta\gamma,\alpha_1\beta\gamma\},\qquad
\{\alpha_1\alpha_2\beta\gamma\}.
\end{gather*}
The last entry is omitted since the final entry in the bottom row is $0$. Note that indeed all of the above forms are pullbacks of harmonic forms from~$\tilde{M}$ up to some correction terms (the harmonic form is always the first element of the sum) which appear in the first two forms in the fourth line. In fact, these are precisely those harmonic forms from $\tilde{M}$ which are not excluded by their cohomology classes falling into the image of ${\rm d}_2$ and ${\rm d}_3$. These forms in turn are elementary to compute from the computation of the cohomology of~$\tilde{M}$ using our spectral sequence since there the lower row again consists of those harmonic forms from $\mathbb{T}^4$ which where not excluded by the image of ${\rm d}_2$ while the top row is just the Hodge star in~$\tilde{M}$ of these elements. The correction terms make it so that the form remains $\delta$-closed and can themselves be computed by a staircase argument applied to the Hodge star of the given form. In these cases they appear due to the fact that the Hodge star of the given basic harmonic forms was not closed with respect to ${\rm d}^{(2,-1)}$.

The harmonic forms of~$M$ corresponding to classes induced by the top row are easily computed by applying the Hodge star to the above elements. Hence, up to sign, we have
\begin{gather*}
\begin{split}
&\{\eta_1\eta_2\sigma\},\\
&
 \{\eta_1\eta_2\alpha_1\beta-\eta_1\gamma\alpha_1\alpha_2,\eta_1\eta_2\alpha_2\beta-\eta_2\gamma\alpha_1\alpha_2,\eta_1\eta_2\alpha_1\alpha_2,\eta_1\eta_2\sigma\gamma,\eta_1\eta_2\alpha_1\sigma,\eta_1\eta_2\alpha_2\sigma\},\\
&\{\eta_1\eta_2\alpha_2\sigma\gamma,\eta_1\eta_2\alpha_1\sigma\gamma,\eta_1\eta_2\beta\sigma\gamma,\eta_1\eta_2\alpha_1\alpha_2\beta,\eta_1\eta_2\alpha_1\alpha_2\sigma\},\\
&\{\eta_1\eta_2\beta\alpha_2\sigma\gamma,\eta_1\eta_2\alpha_1\alpha_2\sigma\gamma, \eta_1\eta_2\alpha_1\beta\sigma\gamma,\eta_1\eta_2\alpha_1\alpha_2\beta\gamma\},\qquad
\{\eta_1\eta_2\beta\alpha_1\alpha_2\sigma\gamma\}.
\end{split}
\end{gather*}
This time we omit the first entry in the row which is~$0$ (consequently, the first listed set of vectors corresponds to the second entry in the top row).

 To compute the harmonic forms corresponding to the middle row, it suffices to compute the forms from the second and third entry in that row (the remaining entries can be computed via the Hodge star or vanish). The harmonic forms corresponding to the second term are
\[\{\eta_1\alpha_1,\eta_2\alpha_2,\eta_1\sigma,\eta_2\sigma,\eta_1\beta-\gamma\alpha_1, \eta_2\beta-\gamma\alpha_2,\eta_1\alpha_2+\eta_2\alpha_1\}.\]
The list of harmonic forms for the third entry is quite long and hence for the sake of clarity we use the convention that replacement of the index next to $\eta$ or $\alpha$ by~$i$ means that the form is an~element of this set for both values of~$i$. Keeping this convention in mind, we get
\[\{\eta_i\sigma\alpha_i,\eta_i\beta\alpha_i,\eta_i\alpha_1\alpha_2, \eta_i\sigma\gamma,\eta_1\alpha_2\beta+\gamma\alpha_1\alpha_2,\eta_2\alpha_1\beta-\gamma\alpha_1\alpha_2, \eta_2\sigma\alpha_1+\eta_1\sigma\alpha_2\}.\]
The final three elements of this set require some explaining. The last one simply comes from the fact that both $\eta_2\sigma\alpha_1$ and $\eta_1\sigma\alpha_2$ are closed while $\eta_2\sigma\alpha_1-\eta_1\sigma\alpha_2$ is in the image of~${\rm d}$. Since harmonic forms are perpendicular to the image of ${\rm d}$ (by the Hodge decomposition), we get the above form as a result. The remaining two problematic entries simply have correction terms stemming from the fact that even though the first components of this sum induce ${\rm d}_2$-closed elements on $E^{p,q}_2$ their ${\rm d}^{(2,-1)}$ image is (up to sign) $\sigma\alpha_1\alpha_2\beta$. Hence, the correction term $\pm\gamma\alpha_1\alpha_2$ needs to be added to negate that.

 Finally, the sets corresponding to the remaining nonzero entries are just the Hodge stars of the previous two sets. Hence, we have
\begin{gather*}
\{\eta_i\sigma\gamma\alpha_i,\eta_i\beta\gamma\alpha_i,\eta_i\gamma\beta\sigma,\eta_i\alpha_1\alpha_2\beta, \eta_2\alpha_1\sigma\gamma+\eta_1\eta_2\beta\sigma,\eta_1\alpha_2\sigma\gamma-\eta_1\eta_2\beta\sigma,\eta_2\beta\gamma\alpha_1+\eta_1\beta\gamma\alpha_2\},\\
\{\eta_1\alpha_1\gamma\beta\sigma,\eta_2\alpha_2\gamma\beta\sigma,\eta_1\gamma\beta\alpha_1\alpha_2,\eta_2\gamma\beta\alpha_1\alpha_2, \eta_1\gamma\sigma\alpha_1\alpha_2+\eta_1\eta_2\alpha_1\sigma\beta,\\
\qquad{}\eta_2\gamma\sigma\alpha_1\alpha_2+\eta_1\eta_2\alpha_2\sigma\beta, \eta_1\sigma\alpha_2\gamma\beta+\eta_2\sigma\alpha_1\gamma\beta\}.
\end{gather*}
\end{egg}
The above example allows us to see both a pattern which further supports the procedure as well as possible technical difficulties when trying to employ it. The supporting pattern comes down to the observation that the basic part of the leading terms (not the correction terms) in each of the harmonic forms above are the basic harmonic forms from the appropriate classes. On the other hand, the example above already demonstrates that harmonic forms corresponding to elements of $E^{p,q}_2$ can have other non-zero parts with respect to the splitting and that these other parts will in general not be harmonic. The above example along with the preceding discussion and the classification of basic harmonic forms for $\mathcal{S}$-structures (cf.~\cite{My2}) suggests the existence of a link between basic harmonic forms and harmonic forms for locally free isometric abelian Lie group actions provided by our spectral sequence. We believe this line of thought to be worth further investigation. For this purpose, we pose the following question the answer to which should clarify how general this link is.

\begin{ques} We can use the bigrading given by the splitting $TM=T\mathcal{F}\oplus T\mathcal{F}^{\perp}$ to fix an~isomorphism from $E^{p,q}_{\infty}$ to $H^{\bullet}_{\rm dR}(M)$ by extending the given $(p,q)$-form uniquely by elements in $\operatorname{Ker}(\delta_b)$ downward along the diagonal. Under what additional conditions is the $(p,q)$-part of the harmonic representative of the image (through this isomorphism) of an element from $E^{p,q}_{\infty}$ basic harmonic?
\end{ques}

Possible uses of this approach include finding examples of geometrically formal manifolds (cf.~\cite{H,K}).

\section{Relevance to the spectral sequence for Riemannian foliations}
In this section, we present how our main result extends to the spectral sequence of the Riemannian foliation $\mathcal{F}$ (by orbits of the action). It is based on the following observation conveyed to us by the anonymous referees, to whom we are sincerely grateful.

\begin{tw}\label{Comp} Under the assumptions of Theorem~$\ref{Main}$, $E^{p,q}_2$ is isomorphic to the second page of the spectral sequence of the foliation $\mathcal{F}$ $($given by the orbits of the action$)$.
\end{tw}

\begin{proof} Let $i\colon \Omega^{\bullet}_G(M)\to\Omega^{\bullet}(M)$ and $\rho\colon \Omega^{\bullet}(M)\to\Omega^{\bullet}_G(M)$ be the inclusion and the averaging over the action of $\overline{G}$, respectively. It is immediate that $\rho\circ i={\rm Id}$ so in particular the same is true for the map induced on the second page.

Hence, it suffices to show that $i\circ\rho$ induces the identity on the second page as well. This can be done by giving a precise chain homotopy between $(i\circ\rho)^*$ and ${\rm Id}$ decreasing the filtration degree by $1$ which can be constructed as follows. Consider a star-shaped fundamental domain~$F$ centered at $0$ in the universal cover $\tilde{G}\cong\mathbb{R}^{N}$ of $\overline{G}$ along with the lift $\tilde{\mu}$ of the normalized Haar measure on $\overline{G}$. Consider the action $T\colon \tilde{G}\times M\to M$ lifted from the action of $\overline{G}$ and assign to every element $a\in F$ the function $H_a\colon M\times I\to M$ defined by $H_a(x,t)=T(ta,x)$. This family varies smoothly with $a\in F$. The chain homotopy can be given by the formula
\[K\omega=\int_F\biggl(\int_0^1H^*_a\omega\biggr){\rm d}\tilde{\mu}(a),\]
where the inner integral denotes the integration over the fiber $[0,1]$. This is obviously a chain homotopy between the desired maps since it is the composition of the standard chain homotopy between the maps induced by the inclusions $i_0,i_1\colon M\to M\times I$ as $M\times\{0\}$ and $M\times\{1\}$ respectively with a chain map $\Omega^{\bullet}(M)\to\Omega^{\bullet}(M\times I)$ whose restrictions to $0$ and $1$ are the desired maps. The proof is finished by noting that this chain homotopy decreases the filtration degree by at most $1$ (due to integration over $[0,1]$).
\end{proof}

Consequently, all the subsequent pages are isomorphic as well and hence our main result immediately extends to this spectral sequence.

\begin{cor}Let $(M^{n+s},g)$ be a compact Riemannian manifold with a locally free isometric action of a connected abelian Lie group $G$ $($of dimension $s)$ such that $\{{\rm d}\eta_1,\dots,{\rm d}\eta_s\}$ has basic cohomological rank $k$ $($i.e., $\dim \langle [{\rm d}\eta_1], \dots, [{\rm d}\eta_s] \rangle = k$ in $H^2(M\slash\mathcal{F}))$. Then, the spectral sequence of the foliation $\mathcal{F}$ $($given by the orbits of the action$)$ degenerates at the latest at the ${(k+2)}$-th page.\looseness=-1
\end{cor}

\begin{rem} We note that Theorem~\ref{Comp} can be generalized to the case when $G$ is not abelian (stating that the second page of the spectral sequence of invariant forms is isomorphic to the second page of the spectral sequence of the foliation given by the orbits of the action) by working with the Lie algebra of $\overline{G}$ in place of the universal covering and lifting $\mu$ through the exponential map. This however is inconsequential to the results presented here.
\end{rem}

\subsection*{Acknowledgments} We express our gratitude to the anonymous referees for pointing out a link between the spectral sequence of invariant forms and the classical spectral sequence of a~Riemannian foliation. This link led to the contents of the final section of the article and significantly increases the scope of our results.


\pdfbookmark[1]{References}{ref}
\LastPageEnding

\end{document}